\documentclass[11pt]{article}
\usepackage{a4wide}

\begin{document}
\newcommand{\be}{\begin{equation}}
\newcommand{\ee}{\end{equation}}
\newcommand{\la}{\label}
\newcommand{\bea}{\begin{eqnarray}}
\newcommand{\eea}{\end{eqnarray}}
\newcommand{\n}{\nonumber}
\newcommand{\nn}{\nonumber \\}
\setcounter{page}0

\title{On a difference equation for generalizations of Charlier
polynomials}
\author{Herman Bavinck \and Roelof Koekoek}
\date{ }
\maketitle

\begin{abstract}
In this paper we obtain a set of polynomials which are orthogonal with respect
to the classical discrete weight function of the Charlier polynomials at
which an extra point mass at $x=0$ is added. We construct a difference
operator of infinite order for which these new discrete orthogonal
polynomials are eigenfunctions.

\vspace{5mm}


\vspace{5mm}

\end{abstract}

\vspace{2cm}

\begin{center}
\begin{tabular}{l}
Delft University of Technology\\
Faculty of Technical Mathematics and Informatics\\
P.O. Box 5031\\
2600 GA Delft\\
The Netherlands\\
Tel : (31-15) 785822 / 787218\\
Fax : (31-15) 787245\\
E-mail : bavinck@twi.tudelft.nl / koekoek@twi.tudelft.nl
\end{tabular}
\end{center}

\newpage

\section{Introduction}

In \cite{Lagdv} J. Koekoek and R. Koekoek found a differential equation of
the form
$$N\sum_{i=0}^{\infty}a_i(x)y^{(i)}(x)+xy''(x)+(\alpha+1-x)y'(x)+ny(x)=0$$
for the polynomials $\left\{L_n^{\alpha,N}(x)\right\}_{n=0}^{\infty}$ which
are orthogonal on the interval $[0,\infty)$ with respect to the weight function
$$\frac{1}{\Gamma(\alpha+1)}x^{\alpha}e^{-x}+N\delta(x),\;\alpha>-1,
\;N\ge 0.$$
The coefficients $\left\{a_i(x)\right\}_{i=1}^{\infty}$ are independent of
the degree $n$ and $a_0(x)$ is independent of $x$.
When $N>0$ this differential equation is of infinite order in general and for
nonnegative integer values of the parameter $\alpha$ the order reduces to
$2\alpha+4$.

In \cite{sym} R. Koekoek also found a similar differential equation for the
symmetric generalized ultraspherical polynomials
$\left\{P_n^{\alpha,\alpha,M,M}(x)\right\}_{n=0}^{\infty}$ which are
orthogonal on the interval $[-1,1]$ with respect to the weight function
$$\frac{\Gamma(2\alpha+2)}{2^{2\alpha+1}\left\{\Gamma(\alpha+1)\right\}^2}
(1-x^2)^{\alpha}+M\left[\delta(x+1)+\delta(x-1)\right],\;\alpha>-1,
\;M\ge 0.$$

For more details concerning these generalized Jacobi polynomials
$\left\{P_n^{\alpha,\beta,M,N}(x)\right\}_{n=0}^{\infty}$ the reader is
referred to \cite{Koorn}.

In \cite{Erice} R.A. Askey posed the problem of finding difference equations
of a similar form for generalizations of discrete orthogonal polynomials
which are orthogonal with respect to the classical weight function together
with an extra point mass at the point $x=0$.

In this paper we solve this problem for generalizations of the classical
Charlier polynomials.

In fact, we look for difference equations of the form
\be\la{DV}N\sum_{i=0}^{\infty}A_i(x)\Delta^iy(x)+
x\Delta\nabla y(x)+(a-x)\Delta y(x)+ny(x)=0\ee
satisfied by the polynomials $\left\{C_n^{a,N}(x)\right\}_{n=0}^{\infty}$
which are orthogonal with respect to the inner product
$$<f,g>=\sum_{x=0}^{\infty}\frac{e^{-a}a^x}{x!}f(x)g(x)+Nf(0)g(0),
\;a>0,\;N\ge 0.$$
In this paper we give a constructive method to obtain the coefficients
$\left\{A_i(x)\right\}_{i=0}^{\infty}$ in the difference equation (\ref{DV})
and we will show that if $N>0$ the order of this difference equation turns
out to be infinite for all values of the parameter $a>0$.

In \cite{BavHaer} the similar problem for generalizations of the Meixner
polynomials is treated.

In \cite{Bav} H. Bavinck introduced Sobolev-type generalizations of the
Charlier polynomials and in \cite{Bav2} it is shown that these are
eigenfunctions of a difference operator of infinite order as well.

\section{Definitions and notations}

We will use the following definition of the classical Charlier polynomials
\be\la{def}C_n^{(a)}(x):=\sum_{k=0}^n{x \choose k}\frac{(-a)^{n-k}}{(n-k)!}
=\frac{(-a)^n}{n!}{}_2F_0\left(\left.\begin{array}{c}-n,-x\\-\end{array}
\right|-\frac{1}{a}\right),\;n=0,1,2,\ldots.\ee
This definition is slightly different from the one given in \cite{Chihara}
but this one turns out to be very convenient in this work. For further
details concerning the classical Charlier polynomials the reader is referred
to \cite{Chihara} anyway.

We remark that we have
\be\la{Lag}C_n^{(a)}(x)=L_n^{(x-n)}(a),\;n=0,1,2,\ldots\ee
where $L_n^{(\alpha)}(x)$ denotes the Laguerre polynomial defined by
$$L_n^{(\alpha)}(x):=\frac{1}{n!}\sum_{k=0}^n(-n)_k(\alpha+k+1)_{n-k}
\frac{x^k}{k!},\;n=0,1,2,\ldots.$$

Further we define the following difference operators~:
\be\la{delta}\Delta f(x):=f(x+1)-f(x)\ee
and
\be\la{nabla}\nabla f(x):=f(x)-f(x-1).\ee

Then we find
\be\la{diff}\Delta C_n^{(a)}(x)=C_{n-1}^{(a)}(x),\;n=1,2,3,\dots.\ee

We also have
$$C_n^{(a)}(0)=\frac{(-a)^n}{n!},\;n=0,1,2,\ldots$$
and
$$C_n^{(a)}(-1)=(-1)^ne_n^a,\;n=0,1,2,\ldots\;\mbox{ where }\;
e_n^a:=\sum_{k=0}^n\frac{a^k}{k!}.$$

The classical Charlier polynomials are discrete orthogonal polynomials which
satisfy the orthogonality relation given by
$$\sum_{x=0}^{\infty}\frac{e^{-a}a^x}{x!}C_m^{(a)}(x)C_n^{(a)}(x)=
\frac{a^n}{n!}\delta_{mn},\;a>0,\;m,n=0,1,2,\ldots,$$
where $\delta_{mn}$ denotes the Kronecker delta.

They satisfy a second order difference equation which can be written in the
form
$$x\Delta\nabla y(x)+(a-x)\Delta y(x)+ny(x)=0,\;y(x)=C_n^{(a)}(x).$$
By using the definition of the difference operators (\ref{delta}) and
(\ref{nabla}) we may rewrite this as
\be\la{dv}ay(x+1)+(n-a-x)y(x)+xy(x-1)=0,\;y(x)=C_n^{(a)}(x).\ee

From the generating function
$$e^{-at}(1+t)^x=\sum_{n=0}^{\infty}C_n^{(a)}(x)t^n$$
we easily obtain
$$\sum_{m=0}^{\infty}C_m^{(a)}(x)t^m\sum_{j=0}^{\infty}C_j^{(-a)}(-x)t^j
=e^{-at}(1+t)^xe^{at}(1+t)^{-x}=1.$$
Hence
$$\sum_{m=0}^kC_m^{(a)}(x)C_{k-m}^{(-a)}(-x)=\left\{\begin{array}{ll}1,&k=0\\
\\0,&k=1,2,3,\ldots\end{array}\right.$$
This can also be written as
\be\la{cru}\sum_{k=j}^iC_{i-k}^{(a)}(x)C_{k-j}^{(-a)}(-x)=\delta_{ij},
\;j\le i,\;i,j=0,1,2,\ldots.\ee

Formula (\ref{cru}) plays an important role in the sections 4 and 5 of this
paper.

We will use another elegant formula which can be obtained from the
generating function. We have for arbitrary real $p$
$$\sum_{n=0}^{\infty}C_n^{(a)}(x+p)t^n=e^{-at}(1+t)^{x+p}=
(1+t)^p\sum_{m=0}^{\infty}C_m^{(a)}(x)t^m.$$
Hence
$$C_n^{(a)}(x+p)=\sum_{k=0}^n{p \choose k}C_{n-k}^{(a)}(x),
\;n=0,1,2,\ldots.$$
The special case $p=-1$ reads
\be\la{som}C_n^{(a)}(x-1)=\sum_{k=0}^n(-1)^kC_{n-k}^{(a)}(x),
\;n=0,1,2,\ldots.\ee
The special case $p=n$ can be written as
\be\la{hulp}C_n^{(a)}(x+n)=\sum_{k=0}^n{n \choose k}C_k^{(a)}(x),
\;n=0,1,2,\ldots.\ee

\section{Generalizations of the Charlier polynomials}

Let ${\bf P}$ denote the space of all real polynomials with real
coefficients. In this section we will determine a set of polynomials
which are orthogonal with respect to the inner product
\be\la{ip}<f,g>=\sum_{x=0}^{\infty}\frac{e^{-a}a^x}{x!}f(x)g(x)+Nf(0)g(0),
\;a>0,\;N>0\;\mbox{ and }\;f,g\in{\bf P}.\ee

If we denote this set of polynomials by
$\left\{C_n^{a,N}(x)\right\}_{n=0}^{\infty}$, where
degree$[C_n^{a,N}(x)]=n$, then we will show that
coefficients $A_n$ and $B_n$ can be determined in such a way that these
polynomials can be written in the form
$$C_n^{a,N}(x)=A_nC_n^{(a)}(x)+B_nC_n^{(a)}(x-1).$$

Suppose that $n\ge 2$ and
$$p(x)=xq(x) \;\mbox{ with degree}[q(x)]\le n-2.$$
Then we easily obtain by using the orthogonality property of the classical
Charlier polynomials
$$<p(x),C_n^{a,N}(x)>=B_n\sum_{x=0}^{\infty}\frac{e^{-a}a^x}{x!}
xq(x)C_n^{(a)}(x-1)
=aB_n\sum_{x=0}^{\infty}\frac{e^{-a}a^x}{x!}q(x+1)C_n^{(a)}(x)=0.$$
Hence, $A_n$ and $B_n$ must satisfy for $n\ge 1$
$$0=<1,C_n^{a,N}(x)>=B_n\sum_{x=0}^{\infty}\frac{e^{-a}a^x}{x!}
C_n^{(a)}(x-1)+NA_nC_n^{(a)}(0)+NB_nC_n^{(a)}(-1).$$

Now we use (\ref{som}) and the orthogonality property of the classical
Charlier polynomials to obtain
$$\sum_{x=0}^{\infty}\frac{e^{-a}a^x}{x!}C_n^{(a)}(x-1)=(-1)^n.$$
Hence
$$NA_nC_n^{(a)}(0)+\left[(-1)^n+NC_n^{(a)}(-1)\right]B_n=0,
\;n=1,2,3,\ldots.$$
So we may choose
$$A_n=1+N(-1)^nC_n^{(a)}(-1)\;\mbox{ and }\;
B_n=N(-1)^{n-1}C_n^{(a)}(0),\;n=0,1,2,\ldots$$
which leads to the following proposition.

\vspace{3mm}

{\bf Proposition.} {\em The generalized Charlier polynomials
$\left\{C_n^{a,N}(x)\right\}_{n=0}^{\infty}$ which are orthogonal with
respect to the inner product (\ref{ip}) can be defined by
\be\la{Def}C_n^{a,N}(x)=\left[1+N(-1)^nC_n^{(a)}(-1)\right]C_n^{(a)}(x)
-N(-1)^nC_n^{(a)}(0)C_n^{(a)}(x-1),\;n=0,1,2,\ldots.\ee}

\vspace{3mm}

Note that we have chosen $C_0^{a,N}(x)=C_0^{(a)}(x)=1$. Further we remark
that we can write
$$C_n^{a,N}(x)=\left[1+N(-1)^{n-1}C_{n-1}^{(a)}(-1)\right]C_n^{(a)}(x)
+N(-1)^nC_n^{(a)}(0)\Delta C_n^{(a)}(x-1),$$
for $n\ge 1$ since we have
$$C_n^{(a)}(0)-C_n^{(a)}(-1)=C_{n-1}^{(a)}(-1),\;n=1,2,3,\ldots.$$

For convenience we define $C_{-1}^{(a)}(x)\equiv 0$ in the sequel. Now the
latter definition holds for all $n\in\{0,1,2,\ldots\}$.

\section{The difference equation}

We try to find a difference equation of the form (\ref{DV}) for the
polynomials $\left\{C_n^{a,N}(x)\right\}_{n=0}^{\infty}$ found in
the preceding section and given by (\ref{Def}), where the coefficients
$\left\{A_i(x)\right\}_{i=1}^{\infty}$ are arbitrary functions of $x$
independent of the degree $n$. Since we want the polynomials
$\left\{C_n^{a,N}(x)\right\}_{n=0}^{\infty}$ to be eigenfunctions of a
difference operator we assume that $A_0(x)$ does not depend on $x$.

So we set
$$y(x)=C_n^{a,N}(x)=\left[1+N(-1)^nC_n^{(a)}(-1)\right]C_n^{(a)}(x)
-N(-1)^nC_n^{(a)}(0)C_n^{(a)}(x-1)$$
and substitute this in the difference equation (\ref{DV}).
Then we find by means of the form (\ref{dv}) of the difference equation for
the classical Charlier polynomials
\bea & &N\left[1+N(-1)^nC_n^{(a)}(-1)\right]
\sum_{i=0}^{\infty}A_i(x)\Delta^iC_n^{(a)}(x)+{}\nn
& &{}\hspace{1cm}{}-N^2(-1)^nC_n^{(a)}(0)
\sum_{i=0}^{\infty}A_i(x)\Delta^iC_n^{(a)}(x-1)+{}\nn
& &{}\hspace{1cm}{}-N(-1)^nC_n^{(a)}(0)
\left[aC_n^{(a)}(x)+(n-a-x)C_n^{(a)}(x-1)+xC_n^{(a)}(x-2)\right]=0.
\n\eea

By using (\ref{delta}), (\ref{diff}) and the difference equation (\ref{dv})
we obtain
$$aC_n^{(a)}(x)+(n-a-x)C_n^{(a)}(x-1)+xC_n^{(a)}(x-2)=-C_{n-1}^{(a)}(x-2).$$
Hence
\bea & &N\left[1+N(-1)^nC_n^{(a)}(-1)\right]
\sum_{i=0}^{\infty}A_i(x)\Delta^iC_n^{(a)}(x)+{}\nn
& &{}\hspace{1cm}{}-N^2(-1)^nC_n^{(a)}(0)
\sum_{i=0}^{\infty}A_i(x)\Delta^iC_n^{(a)}(x-1)+
N(-1)^nC_n^{(a)}(0)C_{n-1}^{(a)}(x-2)=0.\n\eea

This formula must be valid for all values of $a>0$ and $N>0$. The left-hand
side is a polynomial in $N$. So each coefficient of this polynomial has to
be zero. This implies that
\be\la{form0}C_n^{(a)}(-1)\sum_{i=0}^{\infty}A_i(x)\Delta^iC_n^{(a)}(x)-
C_n^{(a)}(0)\sum_{i=0}^{\infty}A_i(x)\Delta^iC_n^{(a)}(x-1)=0\ee
and
$$\sum_{i=0}^{\infty}A_i(x)\Delta^iC_n^{(a)}(x)+
(-1)^nC_n^{(a)}(0)C_{n-1}^{(a)}(x-2)=0.$$
This can be simplified to
\be\la{form1}\sum_{i=0}^{\infty}A_i(x)\Delta^iC_n^{(a)}(x)
=(-1)^{n-1}C_n^{(a)}(0)C_{n-1}^{(a)}(x-2)\ee
and
\be\la{form2}\sum_{i=0}^{\infty}A_i(x)\Delta^iC_n^{(a)}(x-1)
=(-1)^{n-1}C_n^{(a)}(-1)C_{n-1}^{(a)}(x-2),\ee
since $C_n^{(a)}(0)\ne 0$ and $C_n^{(a)}(-1)\ne 0$.

We will show that (\ref{form1}) and (\ref{form2}) have a unique solution for
the coefficients $\left\{A_i(x)\right\}_{i=0}^{\infty}$ which gives rise to
the following theorem.

\vspace{3mm}

{\bf Theorem 1.} {\em The generalized Charlier polynomials
$\left\{C_n^{a,N}(x)\right\}_{n=0}^{\infty}$ satisfy a unique difference
equation of the form (\ref{DV}) where
\be\la{Anul}A_0(x):=A_0(n,a)=(-1)^{n-1}C_{n-1}^{(a)}(-2),\;n=0,1,2,\ldots\ee
and
\bea\la{Ai}& &A_i(x):=A_i(a,x)\nn
&=&\sum_{k=1}^i(-1)^kC_{i-k}^{(-a)}(-x+1)\times{}\nn
& &{}\hspace{1cm}{}\times
\left[C_k^{(a)}(-1)C_k^{(a)}(x-2)-C_k^{(a)}(-2)C_k^{(a)}(x-1)\right],\;
i=1,2,3,\ldots.\eea}

\vspace{3mm}

The proof of this theorem can be found in the next section.
Here the formula (\ref{cru}) is important. Formula (\ref{cru}) can be stated
in other words as follows. If we define the matrix $T:=(t_{ij})_{i,j=1}^n$
($n\ge 1$) with entries
$$t_{ij}:=\left\{\begin{array}{ll}C_{i-j}^{(a)}(x), & j\le i\\0, & j > i
\end{array}\right.$$
then this matrix $T$ is a triangular matrix with determinant 1 and the
inverse $U$ of this matrix is given by $T^{-1}:=U=(u_{ij})_{i,j=1}^n$ with
entries
$$u_{ij}:=\left\{\begin{array}{ll}C_{i-j}^{(-a)}(-x), & j\le i\\0, & j > i.
\end{array}\right.$$

The difference equation given by (\ref{DV}), (\ref{Anul}) and (\ref{Ai}) is
of infinite order for all values of the parameter $a$. This can be seen as
follows. From (\ref{Ai}) it is clear that degree$\left[A_i(x)\right]\le i$
for all $i=1,2,3,\ldots$. Now we compute the coefficient $h_i$ of $x^i$
in the polynomial $A_i(x)$. By using the definition (\ref{def}) we easily see
that
$$C_n^{(a)}(x)=\frac{1}{n!}x^n+{}\;\mbox{ lower order terms.}$$
Hence, from (\ref{Ai}) we find for $i=1,2,3,\ldots$ by using (\ref{hulp})
and (\ref{diff})
\bea h_i&=&\sum_{k=1}^i(-1)^k\frac{(-1)^{i-k}}{(i-k)!}
\left[\frac{C_k^{(a)}(-1)}{k!}-\frac{C_k^{(a)}(-2)}{k!}\right]\nn
&=&\frac{(-1)^i}{i!}\left[\sum_{k=0}^i{i \choose k}C_k^{(a)}(-1)
-\sum_{k=0}^i{i \choose k}C_k^{(a)}(-2)\right]\nn
&=&\frac{(-1)^i}{i!}\left[C_i^{(a)}(i-1)-C_i^{(a)}(i-2)\right]
=\frac{(-1)^i}{i!}C_{i-1}^{(a)}(i-2).\n\eea
This shows that the difference equation given by (\ref{DV}), (\ref{Anul})
and (\ref{Ai}) is of infinite order, since we have by using (\ref{def}) and
(\ref{Lag})
$$C_{i-1}^{(a)}(i-2)=L_{i-1}^{(-1)}(a),\;i=1,2,3,\ldots$$
or
$$C_{i-1}^{(a)}(i-2)=-\frac{a}{i-1}C_{i-2}^{(a)}(i-1)=
-\frac{a}{i-1}L_{i-2}^{(1)}(a),\;i=2,3,4,\ldots.$$
Only when $a$ is a zero of some Laguerre polynomial one of the leading
coefficients $h_i$ might be zero (for some value of $i$), but in that case
we have $h_{i+1}\ne 0$, since two consecutive Laguerre polynomials have
interlacing zeros.

Moreover, the following theorem shows that the infinite order is
unavoidable.

\vspace{3mm}

{\bf Theorem 2.} {\em Every linear difference equation of the form
$$N\sum_{i=0}^{\infty}B_i(x)\Delta^{k_i}\nabla^{i-k_i}y(x)+
x\Delta\nabla y(x)+(a-x)\Delta y(x)+ny(x)=0,\; 0\le k_i\le i$$
satisfied by the polynomials $\left\{C_n^{a,N}(x)\right\}_{n=0}^{\infty}$
has infinite order.}

\vspace{3mm}

{\bf Proof.} From (\ref{delta}), (\ref{nabla}) and (\ref{diff}) we
easily find that
$$\nabla C_n^{(a)}(x)=C_{n-1}^{(a)}(x-1),\;n=0,1,2,\ldots.$$
This implies, in view of (\ref{diff}), that the leading coefficient of
$\Delta^{k_i}\nabla^{i-k_i}C_n^{(a)}(x)$ is equal to that of
$\Delta^iC_n^{(a)}(x)$. This implies, in view of the relations (\ref{form1})
and (\ref{form2}), that the leading coefficient of each $B_i(x)$ equals that
of each corresponding $A_i(x)$. This proves theorem 2.

\vspace{3mm}

Finally we refer to section 6 for more results concerning the coefficients
$\left\{A_i(x)\right\}_{i=1}^{\infty}$ of the difference equation (\ref{DV}).

\section{Proof of theorem 1.}

In this section we will prove that (\ref{form1}) and (\ref{form2}) have a
unique solution for the coefficients $\left\{A_i(x)\right\}_{i=0}^{\infty}$
with $\left\{A_i(x)\right\}_{i=1}^{\infty}$ independent of $n$ and $A_0(x)$
independent of $x$.

Moreover, we will show that this unique solution is given by (\ref{Anul})
and (\ref{Ai}).

Formula (\ref{form0}) can be written as
$$\sum_{i=0}^{\infty}A_i(x)\left[C_n^{(a)}(-1)\Delta^i C_n^{(a)}(x)-
C_n^{(a)}(0)\Delta^i C_n^{(a)}(x-1)\right]=0.$$
Hence
\bea & &\sum_{i=1}^{\infty}A_i(a,x)\left[C_n^{(a)}(-1)\Delta^i C_n^{(a)}(x)-
C_n^{(a)}(0)\Delta^i C_n^{(a)}(x-1)\right]\nn
&=&A_0(n,a)\left[C_n^{(a)}(0)C_n^{(a)}(x-1)-
C_n^{(a)}(-1)C_n^{(a)}(x)\right].\n\eea
The right-hand side vanishes for $x=0$ and since this must be valid for all
values of $n$ and $a>0$ we conclude step by step that
$A_i(0)=A_i(a,0)=0$ for all $i=1,2,3,\ldots$.

Now (\ref{Anul}) easily follows from (\ref{form1}) or (\ref{form2}), since
$C_n^{(a)}(0)\ne 0$ and $C_n^{(a)}(-1)\ne 0$.

By using (\ref{form1}), (\ref{form2}), (\ref{Anul}) and (\ref{diff}) we
obtain
\bea\sum_{i=1}^nA_i(a,x)C_{n-i}^{(a)}(x-1)&=&
(-1)^{n-1}C_n^{(a)}(-1)C_{n-1}^{(a)}(x-2)-
(-1)^{n-1}C_n^{(a)}(x-1)C_{n-1}^{(a)}(-2)\nn
&=&(-1)^n\left[C_n^{(a)}(-1)C_n^{(a)}(x-2)-
C_n^{(a)}(-2)C_n^{(a)}(x-1)\right].\n\eea
Now we use formula (\ref{cru}) to obtain (\ref{Ai}).

Finally we show that (\ref{form1}) must have the same solution. Since we
have
$$C_n^{(a)}(x)=\Delta C_n^{(a)}(x-1)+C_n^{(a)}(x-1)$$
we find for $n\ge 1$, by using (\ref{diff}), that
\bea\sum_{i=1}^nA_i(x)\Delta^iC_n^{(a)}(x)&=&\sum_{i=1}^nA_i(x)\Delta^{i+1}
C_n^{(a)}(x-1)+\sum_{i=1}^nA_i(x)\Delta^iC_n^{(a)}(x-1)\nn
&=&\sum_{i=1}^{n-1}A_i(x)\Delta^iC_{n-1}^{(a)}(x-1)+\sum_{i=1}^n
A_i(x)\Delta^iC_n^{(a)}(x-1).\n\eea
Since the coefficients $\left\{A_i(x)\right\}_{i=1}^{\infty}$ are
independent of $n$ it is sufficient to show that for $n\ge 1$
\bea & &(-1)^{n-1}\left[C_{n-1}^{(a)}(-1)C_{n-1}^{(a)}(x-2)-
C_{n-1}^{(a)}(-2)C_{n-1}^{(a)}(x-1)\right]+{}\nn
& &{}\hspace{1cm}{}+(-1)^n\left[C_n^{(a)}(-1)C_n^{(a)}(x-2)-
C_n^{(a)}(-2)C_n^{(a)}(x-1)\right]\n\eea
equals
$$(-1)^{n-1}C_n^{(a)}(0)C_{n-1}^{(a)}(x-2)-
(-1)^{n-1}C_n^{(a)}(x)C_{n-1}^{(a)}(-2).$$
The proof of this is straightforward and follows by using
the fact that
$$C_{n-1}^{(a)}(x)=\Delta C_n^{(a)}(x)=C_n^{(a)}(x+1)-C_n^{(a)}(x).$$

\section{Some remarks}

We have proved that the polynomials
$\left\{C_n^{a,N}(x)\right\}_{n=0}^{\infty}$ satisfy a unique difference
equation of the form (\ref{DV}) and that the coefficients
$\left\{A_i(x)\right\}_{i=0}^{\infty}$ are given by (\ref{Anul}) and
(\ref{Ai}). In section 4 we already showed that this difference equation is
of infinite order since
$$A_i(x)=\frac{(-1)^i}{i!}C_{i-1}^{(a)}(i-2)x^i+{}\mbox{ lower order terms}
,\;i=1,2,3,\ldots.$$

Some more details about the coefficients can easily be discovered. For
instance, note that the coefficients are both polynomials in $x$ and in $a$.
As a polynomial in $x$ the coefficient $A_i(x)$ usually has degree $i$.
Moreover, if degree$[A_i(x)]<i$ then we have degree$[A_{i+1}(x)]=i+1$. As a
polynomial in $a$ the coefficient $A_i(x)$ has degree $2i-2$. Moreover, we
have by straightforward calculations
$$A_i(x)=\frac{(-1)^i x}{i!(i-1)!}a^{2i-2}+{}\mbox{ lower order terms},
\;i=1,2,3,\ldots.$$

The classical Charlier polynomials also satisfy a difference equation of
infinite order. This can be shown as follows. Similar to formula (\ref{som})
we have more general
$$y(x-1)=\sum_{i=0}^{\infty}(-1)^i\Delta^i y(x)$$
for polynomials $y(x)$. This implies, by using (\ref{delta}) and
(\ref{dv}), that the classical Charlier polynomials satisfy the infinite
order difference equation given by
$$x\sum_{i=1}^{\infty}(-1)^i\Delta^i y(x)+a\Delta y(x)+ny(x)=0,
\;y(x)=C_n^{(a)}(x).$$

So the difference equation (\ref{DV}) for the generalized Charlier
polynomials $\left\{C_n^{a,N}(x)\right\}_{n=0}^{\infty}$ can also be written
in the form
$$N\sum_{i=0}^{\infty}A_i(x)\Delta^i y(x)+x\sum_{i=1}^{\infty}(-1)^i\Delta^i
y(x)+a\Delta y(x)+ny(x)=0.$$

\noindent
{\Large\bf Acknowledgement.}

\vspace{3mm}

\noindent
We thank H. van Haeringen and the editor A.P. Magnus for their useful
suggestions and comments.

\end{document}